\documentclass{amsart}[12pt]
\usepackage{amssymb,latexsym,amsmath,amsthm}

\ifx\pdfoutput\undefined 
\usepackage{graphicx}
\else
\usepackage[pdftex]{graphicx}
\usepackage{epstopdf}
\fi

\newtheorem{theorem}{Theorem}[section]
\newtheorem{proposition}[theorem]{Proposition}
\newtheorem{lemma}[theorem]{Lemma}

\theoremstyle{definition}

\theoremstyle{remark}

\numberwithin{equation}{section}

\newcommand{\frakp}{{\mathfrak p}}

\newcommand{\bbZ}{{\mathbb Z}}
\newcommand{\bbR}{{\mathbb R}}
\newcommand{\bbQ}{{\mathbb Q}}

\newcommand{\h}{\hbar}
\newcommand{\zbar}{\overline{z}}
\DeclareMathOperator{\tr}{Tr}

\newcommand{\ket}[1]{|#1\rangle}
\newcommand{\bra}[1]{\langle #1 |}

\begin{document}
\openup2pt
\newcommand{\medno}{\medskip\noindent}

\renewcommand{\H}{\hat{H}}
\newcommand{\Hcan}{H_{\text{\tiny can}}}

\bigskip
\bigskip

\title{``Bottom of the well" semi-classical trace invariants}

\author{V. Guillemin}\address{Department of Mathematics\\
Massachusetts Institute of Technology\\ 77 Massachusetts Ave.\ \\
Cambridge, MA 02139.} \email{vwg@mit.edu}
\thanks{V.G. supported in part by NSF grant DMS-0408993.}
\author{T. Paul}\address{D\'epartement de math\'ematiques et applications, \'Ecole  
Normale Sup\'erieure \\ 45 rue d'Ulm - F 75230 Paris cedex 05.}
\email{Thierry.Paul@ens.fr}
\author{A. Uribe}
\address{Department of Mathematics\\
University of Michigan\\ 530 Church Street \\ Ann Arbor, Michigan 48109-1043.}
\email{uribe@umich.edu}
\thanks{A.U. supported in part by NSF grant DMS-0401064.}

\maketitle


\begin{abstract}
Let $\hat H$ be an $\h$-admissible pseudodifferential operator whose principal symbol,
$H$, has a unique non-degenerate global minimum.  We give a simple proof that the
semi-classical asymptotics of the eigenvalues of $\hat H$ corresponding to
the ``bottom of the well"
determine the Birkhoff normal form of $H$ at the minimum.  We treat both the resonant
and the non-resonant cases.
\end{abstract}

\section{Introduction}
Let $X$ be an $n-$dimensional manifold and $\H$ a self-adjoint zeroth order semi-classical
$\Psi$DO acting on the space of half-densities, $|\Omega|^{1/2}(X)$.  We will assume that the 
principal symbol, $H(x,\xi)$, of $\H$ has a unique non-degenerate global minimum, $H=C$, 
at some point $(x_0,\xi_0)$, and that outside a small neighborhood of $(x_0,\xi_0)$ $H$ is bounded from below by $C+\delta$, for some $\delta>0$.  We will also assume that at $(x_0,\xi_0)$ the subprincipal symbol of $\H$ vanishes.  From these assumptions one can deduce that on an interval
\[
C < E < C+\epsilon,\quad \epsilon < \delta,
\]
the spectrum of $\H$ is discrete and consists of eigenvalues:
\begin{equation}\label{a1}
E_i(\h),\quad 1\leq i \leq N(\h),
\end{equation}
where
\begin{equation}\label{a2}
N(\h) \sim (2\pi\h)^{-n}\ \text{Vol }\{\;(x,\xi)\;;\; H(x,\xi)\leq C+\epsilon\;\}.
\end{equation}
In addition, we will make a non-degeneracy assumption on the Hessian of $H(x,\xi)$ at $(x_0,\xi_0)$.
Choose a Darboux coordinate system centered at $(x_0,\xi_0)$ such that
\begin{equation}\label{a2}
H(x,\xi) = C + \sum_{i=1}^n \frac{u_i}{2} (x_i^2 + \xi_i^2) + \cdots .
\end{equation}

In this paper we present a short proof of the following theorem:

\begin{theorem}\label{Main}
Assume that the $u_i$'s are linearly independent over the rationals and that
the subprincipal symbol of $\H$ vanishes at $(x_0,\xi_0)$.
Then the eigenvalues, (\ref{a1}), determine the Taylor series of $H$
at $(x_0,\xi_0)$ up to symplectomorphism or, in other words, determine the {\em Birkhoff canonical form}
of $H$ at $(x_0,\xi_0)$.  
\end{theorem}

Our results are closely related to some recent results of \cite{GU} on the Schr\"odinger
operator, $\H = \h^2\Delta + V$, which show that the ``bottom of the well" spectral asymptotics
determines the Taylor series of $V$ at $x_0$.  They are also related to inverse spectral results of
\cite{Gu}, \cite{Ze}, \cite{GP} and \cite{ESZ}.  In these papers it is shown that if
\begin{equation}\label{a4}
\exp tv_H,\qquad v_H = \sum \frac{\partial H}{\partial \xi_i}\, \frac{\partial\ }{\partial x_i} -
\frac{\partial H}{\partial x_i}\, \frac{\partial\ }{\partial \xi_i}
\end{equation}
is the classical dynamical system on $T^*M$ associated with $H$ and $T_\gamma$ is the period of a periodic
trajectory $\gamma$ of this system, the asymptotic behavior of the wave trace, trace$\Bigl(\exp \frac{-it\H}{\h}\Bigr)$ at
$NT_\gamma$, $N\in\bbZ$, determines the Birkhoff canonical form of (\ref{a4}) in a formal 
neighborhood of $\gamma$.  (In the first two papers we have cited, results of this type are proved 
for standard $\Psi$DOs, and \cite{GP} and \cite{ESZ} are versions of these results in the semiclassical
setting.)   It turns out that if the trajectory, $\gamma$, is replaced by a fixed point of the system (\ref{a4}) and,
in particular, if this fixed point is a non-degenerate minimum of $H(x,\xi)$, the recovery
of the Birkhoff canonical form from the spectral data, (\ref{a1}), can be greatly simplified.  Our goal, in this short note, is to show why.  

We also obtain nearly-optimal results in the resonant case, see \S 4.

\section{The semi-classical Birkhoff canonical form}
We quickly review here the construction of the semi-classical Birkhoff canonical form
of $\H$.  We follow the exposition in \cite{CV} and refer to that paper for details.

Performing a preliminary microlocalization and conjugation by an $\h$-FIO, we can
assume that $\H$ is an operator on $\bbR^n$, and that the global minimum
$(x_0,\xi_0)$ is the origin $(0,0)\in T^*\bbR^n$.  Let us denote by 
\[
[\H ] = \sum_{\alpha, \beta, k} x^\alpha \xi^\beta \h^k
\]
the Taylor series of the full Weyl symbol of $\H$, where the monomial 
$x^\alpha \xi^\beta \h^j$ has degree $|\alpha |+|\beta |+2j$. 
In fact we can assume that
\[
[\H ] = \sum_i \frac{u_i}{2} (x_i^2+\xi_i^2) + \cdots
\]
where the dots indicate terms of degree three and higher.  Notice that
\[
[\H ]|_{\h = 0}  = \text{the Taylor series of the principal symbol of }\H.
\]
Let $\H_2$ denotes the Weyl quantization of 
\[
H_2:=  \sum_i \frac{u_i}{2}(x_i^2+\xi_i^2),
\]
and set
\[
\H = \H_2 + \hat L.
\]
To construct the quantum Birkhoff canonical form of $\H$, one conjugates
$\H$ by suitable Fourier integral operators in order to successively
make higher-order terms in $L$ commute with $\H_2$.  
The resulting series is the quantum Birkhoff canonical form, $\Hcan$
of $H$.

The non-resonance condition implies that we can write
$\Hcan$ in the form:
\begin{equation}\label{a7}
\widehat\Hcan = \hat H_2 + F(P_1,\ldots ,P_n,\h),\quad P_i =
\h^2 D_i^2 + x_i^2
\end{equation}
with $F$ an $\h$-admissible symbol whose Taylor series
is of the form
\begin{equation}\label{a8}
[F]= \sum_{|r|\geq 1} c_r(\h) p^r,
\end{equation}
where $ p_i = \xi_i^2 + x_i^2 $, $r=(r_1,\ldots ,r_n)$,
\begin{equation}\label{I}
c_r(\h) = \h^{|r|-1}\Bigl( c_{r,0} + \cdots\Bigr)
\end{equation}
and $ c_r(0)= 0\ \text{for}\ |r|=1$ (so that all
the monomials in $[F]-H_2$ have degree $\geq 3$).

\bigskip
Theorem \ref{Main} is a direct consequence of the following:

\begin{theorem}\label{Main2}
Under the assumptions of Theorem \ref{Main},
the eigenvalues, (\ref{a1}), determine the semi-classical
Birkhoff canonical form of $\H$.
\end{theorem}

\section{The proof of Theorem \ref{Main2}}

The first step in our argument is more or less identical with that
of \cite{GU}, \cite{GP} and \cite{ESZ}.  Assume without loss of
generality that $C=0$, and let $\rho\in C_0^\infty(\bbR)$ be
equal to one on the interval $[-1/2\,,\,1/2]$ and zero outside the
interval $[-1\,,\,1]$.  Then for $\epsilon$ small the $\rho$-truncated
wave trace
\begin{equation}\label{a5}
\text{trace }\rho\Bigl(\frac{\H}{\epsilon}\Bigr)  \exp \frac{-it\H}{\h}
\end{equation}
is equal modulo $O(\h^\infty)$ to the $\rho$-truncated wave trace
for the Birkhoff canonical form,
\begin{equation}\label{II}
\tr(t,\h) := \text{trace }\rho\Bigl(\frac{\widehat\Hcan}{\epsilon}\Bigr)\
\exp \Bigl(\frac{-it\widehat\Hcan}{\h}\Bigr).
\end{equation}
The truncated wave trace admits an asymptotic expansion 
$\tr (t,\h) \sim a_0(t) + a_1(t)\h +\cdots$ as $\h\to 0$.  This
follows from the method of stationary phase
and fact that for each $t$ the operator 
$\rho(\epsilon^{-1}\hat{H})e^{it\h^{-1}\hat{H}}$ is an
$\h$-Fourier integral operator.  Writing the truncated trace as
an oscillatory integral, for each $t$ the phase has a unique 
critical point, corresponding to the absolute minimum $(x_0,\xi_0)$
which is a fixed point of the classical flow.
Since the cutoff operator
$\rho\Bigl(\frac{\widehat\Hcan}{\epsilon}\Bigr)$
is microlocally equal to the identity in a neighborhood of
$(x_0,\xi_0)$, the asymptotic expansion of the cutoff trace 
is independent of $\rho$, provided $\rho\in C_0^\infty$ is identically
equal to one near zero.

\medskip
The truncated trace of the Birkhoff canonical form equals
\begin{equation}
\tr(t,\h) = \sum_{k\in(\bbZ_{+})^n}
\rho\Bigl(\frac{\Hcan(\h(k+1/2),\h)}{\epsilon}\Bigr)\ 
e^{it\h^{-1} \Hcan(\h(k+1/2),\h)}.
\end{equation}
However, since $\rho$ is identically equal to one in a neighborhood
of zero, as a power series in $\h$
\begin{equation}\label{3.4}
\tr(t,\h) \sim \sum_{k\in(\bbZ_{+})^n}
e^{it\h^{-1} \Hcan (\h(k+1/2),\h)}.
\end{equation}

We now rewrite (\ref{3.4}) in a more amenable form using a variant
of the ``Zelditch trick" (see \cite{Ze}).

\begin{proposition}
For any choice of $\rho$ as above and for $\epsilon$ small,
as $\h\to 0$
\begin{equation}\label{IV}
\tr (t,\h) \sim \sum_{m=0}^\infty \frac{(it)^m}{m!}\ 
\Bigl(\sum_{|r|\geq 1} \h^{|r|-1} c_r(\h)
\Bigl(\frac{1}{t} D_\theta\Bigr)^r \Bigr)^m\ 
\frac{e^{it\frac{1}{2}\sum_j\theta_j}}
{\Pi_j (1-e^{i t\theta_j})}\Bigr|_{\theta = u},
\end{equation}
where 
\[
D_\theta = -i\Bigl(\frac{\partial\ }{\partial\theta_1},\ldots ,
\frac{\partial\ }{\partial\theta_n}\Bigr)
\]
and the right-hand side of (\ref{IV}) is understood as a power series
in $\h$.
\end{proposition}
\begin{proof}
Recalling that $[\hat H_2 , \hat F] = 0$,
\[
\tr(t,\h) \sim 
\sum_{k\in(\bbZ_+)^n}
e^{itu\cdot (k+1/2)}\ 
\bra{k}e^{it\h^{-1}\hat{F}}\ket{k},
\]
where $\{\ket{k}\}$ is an orthonormal basis of eigenvectors of
the canonical $n$-torus representation on $L^2(\bbR^n)$, and
$u\cdot (k+1/2) = \sum_{j=1}^n u_j(k_j+1/2)$.
For each $k$, the Taylor expansion, $[F]$, gives us an 
asymptotic expansion
\begin{equation}
\bra{k}e^{it\h^{-1}\hat{F}}\ket{k} =
\sum_{m=0}^\infty \frac{(it)^m}{\h^m m!} F(\h(k+1/2),\h)^m \sim 
\sum_{m=0}^\infty \frac{(it)^m}{m!}
\Bigl(\sum_r \h^{|r|-1}\,c_r(\h) (k+1/2)^r\Bigr)^m.
\end{equation}
Let us introduce the variables $\theta = (\theta_1,\ldots ,\theta_n)$
and write:
\[
(k+1/2)^r e^{itu\cdot (k+1/2)} = \Bigl(\frac{1}{t} D_\theta\Bigr)^r
e^{it\theta\cdot (k+1/2)} \Bigr|_{\theta = u}.
\]
Then
\[
\tr(t,\h) \sim \sum_{m=0}^\infty
\frac{(it)^m}{m!} 
\sum_{k\in\bbZ_+} \rho\Bigl(\frac{F(\h(k+1/2),\h)}{\epsilon}\Bigr)\
\Bigl(\sum_r \h^{|r|-1}\,
c_r(\h) \Bigl(\frac{1}{t} D_\theta\Bigr)^r\Bigr)^m
e^{it\theta\cdot (k+1/2)} \Bigr|_{\theta = u}.
\]
Finally, for each $m$ (summing a geometric series)
\[
\sum_{k\in\bbZ_+}\Bigl(\sum_r \h^{|r|-1}\,
c_r(\h) \Bigl(\frac{1}{t} D_\theta\Bigr)^r\Bigr)^m
e^{it\theta\cdot (k+1/2)} \Bigr|_{\theta = u} = \qquad\qquad \mbox{ }
\]
\[
\mbox{}\qquad\qquad =
\Bigl(\sum_r \h^{|r|-1}\,
c_r(\h) \Bigl(\frac{1}{t} D_\theta\Bigr)^r\Bigr)^m\ 
\frac{e^{it\frac{1}{2}\sum_j \theta_j}}{\Pi_i(1-e^{i t\theta_j})}
\Bigr|_{\theta = u},
\]
and the result follows.
\end{proof}

\bigskip
We will show that the $m=0$ term in the series on
the right-hand side of (\ref{IV}) suffices to determine
the $c_r(\h)$.  More precisely:

\begin{theorem}\label{One}
From the coefficients of $\h^s$, $s\leq\ell$, in the series
in $\h$
\begin{equation}\label{V}
V(t,\h) = \sum_{|r|\geq 1}
\h^{|r|-1}\, c_r(\h) \Bigl(\frac{1}{t} D_\theta\Bigr)^r \ 
\frac{e^{it\frac{1}{2}\sum_j\theta_j}}
{\Pi_j (1-e^{i t\theta_j})}\Bigr|_{\theta = u}
\end{equation}
one can determine the coefficients of $\h^s$, $s\leq\ell$, in $c_r(\h)$
for all $r$.
\end{theorem}
\begin{proof}
Let $\rho$ be a cutoff function as before, and $\hat\varphi = \rho$.
Integrating (\ref{V}) against $\epsilon^n\varphi(\epsilon t)$
and essentially reversing the previous calculation, we find:
\begin{equation}\label{vt}
\tilde{V}(\epsilon,\h) = 
\sum_k\Bigl(\sum_r \h^{|r|-1} c_r(\h) (k+1/2)^r\Bigr)\ 
\rho\Bigl(\epsilon^{-1}u\cdot (k+1/2)\Bigr) = 
\end{equation}
\[
= \sum_k \h^{-1} F(\h(k+1/2))\
\rho\Bigl(\epsilon^{-1}u\cdot (k+1/2)\Bigr).
\]
Letting
\[
c_r(\h) = \sum_{i=0}^\infty c_{r,i}\, \h^{|r|-1+i},
\]
we can rearrange (\ref{vt}) in increasing powers of $\h$ 
(using the variable $\ell= 2|r|-2+i$ for the exponent of $\h$):
\begin{equation}\label{vt1}
\tilde{V}(\epsilon,\h) =
\sum_{\ell=0}^\infty\h^\ell\
\sum_{j=0}^{[\frac{\ell}{2}]} \sum_{|r|=j+1}\
c_{r,\ell-2j}\ \Bigl(
\sum_k\ 
(k+1/2)^r\rho\Bigl(\epsilon^{-1}u\cdot (k+1/2)\Bigr)\Bigr).
\end{equation}
Now arrange the numbers
\[
u_k = u\cdot (k+1/2),\qquad k\in (\bbZ_{+})^n
\]
in strictly increasing order (which is possible because there
are no repetitions among them):
\begin{equation}
0 < \nu_1 =u_k|_{k =0} < \nu_2 < \cdots .
\end{equation}
Let us write: $\nu_s = u_{k^{(s)}}$.
Now vary $\epsilon$ in (\ref{vt1}), starting with a very small value. 
Gradually increasing $\epsilon$, we can arrange
that the coefficient
of $\h^\ell$ in (\ref{vt1}) is
\[
\sum_{j=0}^{[\frac{\ell}{2}]} \sum_{|r|=j+1}\ 
c_{r,\ell-2j}\ 
\sum_{s=1}^m (k^{(s)}+1/2)^r\ 
\rho\Bigl(\epsilon^{-1}\nu_s\Bigr)
\]
for any given $m$.  Therefore, by an inductive argument on
$m$ we can recover the values of the polynomial
\[
\frakp_\ell (x) = \sum_{j=0}^{[\frac{\ell}{2}]} \sum_{|r|=j+1}\
c_{r,\ell-2j}\ (x+1/2)^r
\]
at all $k\in (\bbZ_{+})^n$.  But these values determine the polynomial,
and therefore its coefficients.
\end{proof}

Now we show that the asymptotic expansion of the trace, $\tr(t,\h)$,
determines $V$:
\begin{theorem}\label{Two}
From the coefficients of $\h^s$, $s\leq\ell$, in the expansion
(\ref{IV}) one can determine the coefficients of $\h^s$, $s\leq\ell$
in the series $V(t,\h)$.
\end{theorem}
\begin{proof}
We proceed by induction on $\ell$.  

The coefficient of $\h$ in $V$ coincides with the coefficient of
$\h$ in (\ref{IV}), since all terms in the sum (\ref{IV}) except the
first are of order $O(\h^m)$, $m>1$.

By theorem \ref{One} the coefficient of $\h$ in $V$ enables us to
determine the coefficient of $\h$ in $c_r(\h)$, and hence the
coefficient of $\h^2$ in the second summand of (\ref{IV}).
But the coefficient of $\h^2$ in the first summand coincides with the
coefficient of $\h^2$ in $V$, so the coefficients of $\h$ and $\h^2$
in (\ref{IV}) determine the coefficient of $\h^2$ in $V$.  It is clear
that this procedure can be continued indefinitely.
\end{proof}

\bigskip
Theorem \ref{Main2} is an immediate consequence of theorems \ref{One}
and \ref{Two}.

\section{The resonant case}
We now consider the case when
the frequencies $u_1,\ldots , u_n$ are not linearly independent
over $\bbQ$.  Following \cite{CV}, let us introduce the number
\[
d = \min\{ |\alpha |, \alpha\in\bbZ^n\setminus\{0\}\; |\;
\alpha\cdot u = 0\}
\]
which is a measure of the rational relations among the frequencies
(here $|\alpha| = \sum_{j=1}^n |\alpha_j|$).  
We will make use below of the following observation:
\begin{lemma}
Among the eigenvalues of $\hat H_2$ of the form: 
\begin{equation}\label{res0}
\lambda_k = k\cdot (u+1/2)\qquad\text{with}
\qquad  |k|<d/2
\end{equation}
there are no repetitions (i.\ e.\ 
the mapping $k\mapsto\lambda_k$ is 1-1
in the range $|k|<d/2$).
\end{lemma}
\begin{proof}
If $k\cdot (u+1/2) = k'\cdot (u+1/2)$, then $(k-k')\cdot u = 0$
and therefore $|k-k'|\geq d$.  The conclusion now follows from the
triangle inequality.
\end{proof}

\bigskip
Continuing to assume that the subprincipal symbol vanishes at the
absolute minimum, the semi-classical
Birkhoff canonical form in the resonant
case has the following structure (see \cite{CV}):
\[
H_{\tiny can} = H_2 + F + K,
\]
where
\begin{enumerate}
\item $F = F(p_1,\ldots ,p_n,\h)$ where $F$ is a polynomial
in all variables of degree at most $[\frac{d-1}{2}]$.
\item $[K]$ is a power series with monomials $\h^j x^\alpha \xi^\beta$
where $|\alpha| + |\beta | + 2j > d$ and $[\hat H_2 , \hat K] = 0$.
\end{enumerate}

In this section we prove the following:
\begin{theorem}\label{Resonant}
If $d$ is even the eigenvalues (\ref{a1}) 
determine the entire semi-classical canonical form, $F(x,\h)$. 
If $d$ is odd, those eigenvalues determine the semi-classical
canonical form except for the monomials of maximal degree,
$[\frac{d-1}{2}]$. 
\end{theorem}

Except for a few additional complications,
the method of proof is the same as in the non-resonant case. 
We begin by checking that the asymptotic expansion of the
truncated trace can be treated by the same methods as before, up to 
a sufficiently high order in $\h$:

\begin{proposition}
In the resonant case, the expansion (\ref{IV}) is valid
modulo $O(\h^{[\frac{d}{2}]})$.
\end{proposition}
\begin{proof}
Once again we write the trace as a sum of diagonal
matrix elements over
a normalized basis $\{\ket{k}\}$ of eigenfunctions of 
the standard representation of the $n$-torus, splitting off the
$H_2$ part (which is possible since $\hat F+\hat K$ commutes with $H_2$): 
\[
\tr(t,\h) \sim
\sum_{k\in(\bbZ_+)^n}
e^{itu\cdot (k+1/2)}\
\bra{k}e^{it\h^{-1}(\hat{F}+\hat{K})}\ket{k}.
\]
We next expand the exponential in its Taylor series. 
We want to show that every term involving $\hat K$
is $O(\h^{[\frac{d}{2}]})$.

A term involving
\[
\bra{k}(\hat F +\hat K)^m\ket{k}
\]
is a sum of terms of the form
\begin{equation}\label{res1}
\bra{k}\hat{F_1}\hat{K_1}\cdots\hat{F_s}\hat{K_s}\ket{k}
\end{equation}
where the $F_j$ are powers of $F$ and the $K_j$ are powers of $K$.
Therefore, the $K_j$ are sums of monomials $\h^j x^\alpha \xi^\beta$
where $2j + |\alpha| + |\beta| > d$, just as is $K$.
Let us express those monomials in terms of raising and lowering operators,
\[
A_{\alpha \beta} = z^\alpha \zbar^\beta,\qquad
z=x+i\xi.
\]
Then (\ref{res1}) is a linear combination of terms of the form
\begin{equation}\label{res2}
\h^{\sum_{i=1}^s j_i}\bra{k}\hat{F_1}\widehat{A_{\alpha^1 \beta^1}}\cdots
\hat{F_s}\widehat{A_{\alpha^s\beta^s}}\ket{k}
\end{equation}
where, for each $i$, 
\[
2j_i+|\alpha^i|+|\beta^i|>d.
\]
Now recall that (i) the $\hat{F}_j$ are diagonal
in the basis $\{\ket{k}\}$ and (ii)
the $\widehat{A_{\alpha \beta} }$
act on the basis vectors by:
\[
\widehat{A_{\alpha \beta}} \ket{k} =
\h^{|\beta|}c_{\alpha\beta}\ket{k+\alpha-\beta}
\]
where $c_{\alpha\beta}$ is a constant whose value we won't need.
Therefore, a diagonal matrix element of the sort
(\ref{res2}) is zero unless 
\[
\sum_{i=1}^s \alpha^i-\beta^i = 0,
\]
in which case (\ref{res2}) is $O(\h^{j+\sum_i |\beta^i|})$
where $j = \sum_i j_i$.
However, $|\alpha^i|+|\beta^i| > d-2j_i$ for each $i$ and
\[
\sum_{i=1}^s \alpha^i-\beta^i = 0\quad \Rightarrow\quad
\Bigl|\sum_{i=1}^s \alpha^i\Bigr| = \Bigl|\sum_{i=1}^s \beta^i\Bigr|.
\]
Therefore, $\sum_i |\beta^i| \geq [sd/2]-j$ and so
(\ref{res2}) is $O(\h^{[sd/2]})$.  It follows that
all diagonal matrix elements to which $\hat{K}$ contributes 
are at least $O(\h^{[\frac{d}{2}]})$.
\end{proof}

\begin{lemma}
$F(\h(x+1/2),\h)$ is a polynomial in $\h$ of degree
at most $[\frac{d-1}{2}]$, and if we write
\[
F(\h(x+1/2),\h) = \sum_{j=0}^{[\frac{d-1}{2}]} h^j F_j(x)
\]
the power series expansion of $\tr (t,\h)$ determines the
values $F_j(k)$ for all $k\in (\bbZ_{+})^n$ such that
$|k|<d/2$, for all $j\leq [\frac{d-1}{2}]$
if $d$ is even and for all $j<[\frac{d-1}{2}]$ if $d$ is odd.
\end{lemma}
\begin{proof}
The first statement follows from the general form of $F$.

By theorems \ref{Two} and \ref{One}, for any $\ell$
the first $\ell$ terms of the expansion of
$\tr (t,\h)$ determine the first $\ell$ terms of 
(\ref{vt}), provided we replace $F$ by $F+K$.
But, by the previous proposition, the 
expansion of $\tr (t, \h)$ mod $O(\h^{[\frac{d}{2}]})$ is 
insensitive to what $K$ is.  Therefore, (\ref{vt})
remains valid mod $O(\h^{[\frac{d}{2}]})$,
where $F$ now stands for 
the part of the canonical form we are determining from
the spectrum. 

If $d$ is even
\[
[\frac{d-1}{2}] < [\frac{d}{2}],
\]
and so
it follows that the expansion of $\tr$ determines the sums
\[
\sum_k \h^{-1} F(\h(k+1/2))\
\rho\Bigl(\epsilon^{-1}u\cdot (k+1/2)\Bigr).
\]
Now we proceed as before, letting $\epsilon$ grow starting at a 
very small value.  Since the eigenvalues (\ref{res0}) are all different,
we can determine the polynomial in $\h$,
$F(\h(x+1/2))$, evaluated at each $k$ with $|k|<d/2$.
If $d$ is odd we must discard the term $F_j$ with $j=[\frac{d-1}{2}]$. 
\end{proof}

\bigskip
Since  $F_j$ is a polynomial of degree at most
$[\frac{d-1}{2}]$, the proof of theorem \ref{Resonant}
is completed by the following result:
\begin{lemma}
Let $f(x_1,\ldots , x_n)$ be a polynomial of degree $N$.  Then
$f$ is completely determined by its values at the points
\[
(k_1+1/2,\ldots , k_n+1/2),
\]
for all $k$ such that $|k|\leq N$ and $k_j\geq 0$.
\end{lemma}
\begin{proof}
The proof is by induction on the number of variables.
The case $n=1$ is trivial.
Assume the result is true for polynomials of $n-1$
variables, and let
\[
f = f_N(x_2,\ldots ,x_n)+ f_{N-1}(x_2,\ldots ,x_n)x_1 +\cdots +
f_0\, x_1^N.
\]
Note that degree $f_i = i$.  

Evaluating $f$ at
$(k+1/2, 1/2,\ldots , 1/2)$, $0\leq k\leq N$ determines 
$f_i(1/2,1/2,\ldots ,1/2)$, $i=0,\ldots N$, and in particular
determines $f_0$.

Evaluating $f-f_0\,x_1^N$ at
$(k+1/2,k_2+1/2,\ldots ,k_n+1/2)$, $0\leq k\leq N-1$,
$k_2+\cdots + k_n\leq 1$ determines
$f_i(k_2+1/2,\ldots ,k_n+1/2)$ at all $k_2+\cdots + k_n\leq 1$ 
and in particular determines $f_1$.

Evaluating $f-f_1\,x_1^{N-1}-f_0\,x^N$ at
$(k+1/2,k_2+1/2,\ldots ,k_n+1/2)$, $0\leq k\leq N-2$,
$k_2+\cdots + k_n\leq 2$ determines
$f_i(k_2+1/2,\ldots ,k_n+1/2)$ at all $k_2+\cdots + k_n\leq 2$
and in particular determines $f_2$.  Etc.

\end{proof}

\bigskip
When $d$ is odd our methods recover the classical
Birkhoff normal form of $H$ except for its monomials
of top degree, $\frac{d-1}{2}$.


\end{document}